\documentclass[12pt]{article}

\usepackage[utf8]{inputenc}
\usepackage[T1]{fontenc}
\usepackage{lmodern}

\usepackage{csquotes} % for \enquote

\usepackage{amsmath, amssymb}
\usepackage{amsthm}

\usepackage{graphicx}
\usepackage{float}
\usepackage{xcolor}

\usepackage[margin=1in]{geometry}

\usepackage[normalem]{ulem}

\usepackage{hyperref}
\hypersetup{
    colorlinks=false,
    linkcolor=blue,
    urlcolor=blue,
    citecolor=blue,
    hidelinks
}

\title{\textbf{Three Approaches to the Problem of Double Discontinuity}}
\author{Manya Raman-Sundström}
\date{\today}

\begin{document}

\maketitle

\section*{Introduction}

In 1872, Felix Klein gave an inaugural speech for his professorship at Erlangen, in which he laid out some of his core values for school mathematics. He said:

\begin{quote}
\small
[W]e, as university teachers, require not only that our students, on completion of their studies, know what is to be taught in the schools. We want the future teacher to stand \emph{above} his subject, that he have a conception of the present state of knowledge in his field, and that he generally be capable of following its further development. Therefore, we hope to lead him far enough that he at least once undertakes an independent research study (Klein, \emph{Antrittsrede}, Rowe tr., 1985, p. 139)\footnote{Original: [W]ir als Universitaetslehrer von unseren Zuhoerern nicht nur, dass sie am Ende Ihrer Studien wissen, was auf der Schule vorzutragen ist. Wir wollen, dass der spaetere Lehrer ueber seinem Stoffe steht, dass er einen Begriff von dem hat, was seine Wissenschaft zur Zeit umfasst, das er in Allgemeinen zu verfolgen vermag, in welchem Sinne sie sich weiter entwickelt.}.
\end{quote}

This was roughly twenty-five years before Klein would receive the commission from the German government to help reform mathematics education across the newly reunited regions of the country\footnote{See Klein 1900.}. And it was thirty years before he, with Rudolf Schimmack, David Hilbert, Wilhelm Lorey, and other colleagues, would draft the Merano syllabus\footnote{Klein and Schimmack, et al, 1905.}, which would become the blueprint for reform in mathematics around the world. Already here, we see the germs of what would become the foundation for Klein's views of teacher education: to give future teachers a higher standpoint.

\subsection*{The problem of ``double discontinuity''}

Klein wanted future teachers to have a higher standpoint to address a fundamental gap between school mathematics and university mathematics. This gap gave rise to two problems, which Klein called the ``double discontinuity'' (Klein, 1908/2016, p.1): students entering university were unprepared for higher education, and after university, they were unprepared to teach school mathematics effectively. Klein wrote:

\begin{quote}
\small
The young university student finds himself, at the outset, confronted with problems, which in no way remind him of the things with which he had been concerned at school. Naturally he forgets all these things quickly and thoroughly. When, after finishing his course of study, he becomes a teacher, he suddenly finds himself expected to teach the traditional elementary mathematics according to school practice; and, since he will be scarcely able, unaided, to discern any connection between this task and his university mathematics, he will soon fall in with the time-honored way of teaching, and his university studies remain only a more or less pleasant memory which has no influence upon his teaching\footnote{``Der junge Universitats-Student findet sich zu Anfang mit Aufgaben konfrontiert, die ihn in keiner Weise an die Dinge erinnern, mit denen er in der Schule zu tun hatte. Selbstverstandlich vergisst er sie rasch und grundlich. Wenn --- nachdem er sein Studium beendet hat --- Lehrer wird, so findet er sich auf einmal in der Lage, die traditionelle Elementarmathematik gemass der schulischen Gepflogenheit zu lehren; und da er ohne fremde Hilfe kaum in der Lage sein wird, eine Verbindung zwischen dieser Aufgabe und seiner Hochschulmathematik zu erkennen, fugt er sich bald der ehrwurdigen Lehrweise, und sein Hochschulstudium bleibt nur eine mehr oder weniger angenehme Erinnerung, die keinen Einfluss auf seinen Unterricht ausubt.'' Original.} (Klein 1908/2016).
\end{quote}

Klein aimed to better align school and university mathematics. The German government required bottom-up reforms\footnote{Report on Schools, 1900.}, so Klein observed classes, spoke with teachers, and developed a course for pre-service teachers. These notes culminated in the three-volume \emph{Elementary Mathematics from a Higher Standpoint} (EMHS), covering the material Klein thought mathematics teachers should know (Klein, 1908/2016). However, the notes do not explicitly solve the double discontinuity problem.

Klein emphasized the teacher's perspective rather than the student's. He wrote:

\begin{quote}
\small
What is required is more interest in mathematics, livelier instruction, and a more spirited treatment of the material!\footnote{''Aber wir verlangen mehr Interesse fur die Mathematik, mehr Leben in ihrem Unterrichte, mehr Geist in ihrer Behandlung!'' (original).} (Klein, \emph{Antrittsrede}, translated Rowe p. 139)
\end{quote}

Thus, Klein described the problem of double discontinuity without prescribing a method for connecting school and university mathematics. The notes provide insight into terms like ``elementary'' and ``higher standpoint'' but no roadmap for bridging the gap\footnote{See Martino et al., 2022, for recent literature on the high school to university transition.}.

This paper reconceptualizes the double discontinuity problem along two dimensions. We examine three approaches—vertical and horizontal—highlighting how horizontal approaches can address the problem in ways vertical ones cannot.

\subsection*{Structure of the paper}

We begin with an excerpt from Klein's \emph{Elementary Mathematics from a Higher Standpoint} on logarithms, exemplifying his \textbf{vertical} development from school to research mathematics\footnote{See Kilpatrick 2019; Allmendinger 2019.}.

Next, we consider \emph{Roots to Research: A Vertical Development of Mathematical Problems} (RR) by Paul and Judith Sally, which also exemplifies vertical development. Their work spans grade-school to research-level mathematics, targeting teachers and serious students.

Finally, we present an \emph{Extended Analysis} of a high school problem, the Box Problem from \emph{Mathematics for High School Teachers---An Advanced Perspective} (MHST) by Dick Stanley et al. This \textbf{horizontal} approach uses high school mathematics but applies deeper reasoning to produce results that are simple, generalizable, and mathematically rich, echoing ideas from Polya (1945) and Boaler (2022).

\vspace{1em}
Detailed descriptions of these examples illuminate the development of mathematical ideas in different contexts, showing how double discontinuity arises and how it might be avoided.

%%put example1 here

\section*{Example 1: Logarithms}

The first example comes from Klein's \emph{Elementarmathematik vom höheren
Standpunkte aus} (EMHS)\footnote{Note that the translation of the title into
English has been problematic. An earlier translation rendered ``höher'' as
``advanced.'' The editor of the updated version, used here, claims the correct
translation is ``higher,'' and argues that the difference is substantial
(Schubring, 2016, p. v).}. The chapter starts with the high school definition
of logarithm, traces its historical development\footnote{Klein's use of history
to teach mathematical concepts was central to his method: ``If you do know
the historical development, your footing will be very insecure'' (Klein 1924).
See Kilpatrick, 2019 for discussion.}, and concludes with a treatment in
complex analysis, which Klein called function theory. Klein begins with
Napier's definition of logarithm as an exponent, moves to Euler's notion of
logarithm as a function, and then discusses modern treatments (for his time)
over positive, integer, real, and complex bases.

Klein's treatment illustrates his pedagogical idea of ``elementarization,''
which is how a complex concept can be distilled to its fundamental nature\footnote{
According to Schubring, Klein's conception of ``elementarization'' was influenced
by D'Alembert's ``éléments des sciences'' in the \emph{Encyclopédie}, which aimed
to distill science to its core and connect it to the whole of science\footnote{D'Alembert,
1755, as quoted in Schubring, 2016.}. This differs from Euclid's \emph{Elements},
used in his time to denote the fundamentals. The term was also used by Apollonius
in his \emph{Conics} (Fried, 2001), but not in exactly the same sense.
}. To elementarize is not to simplify, though it may appear so in a vertical
treatment. An elementary idea might represent an essence that emerges fully
only after considering the concept in all its generality. This treatment
also reflects Klein's idea of \emph{hysteresis}, which refers to letting
concepts mature before introducing them into the school curriculum\footnote{See
Schubring, 2016.}. Once matured, these concepts form the backbone of school
mathematics; ``functions'' is another such concept.

While Klein comments on school-level teaching (for example, he preferred
defining logarithms via integrals rather than as inverses of exponents),
his main goal is to show the coherence of the notion from school to research
mathematics. He writes:

\begin{quote}
\small
I shall by no means address myself to beginners, but I shall take for
granted that you are all acquainted with the main features of the most
important disciplines of mathematics. <\dots> My main task will always
be to show you the mutual connection between problems in various disciplines;
these connections used not to be sufficiently considered in the specialized
lecture courses, and I want more especially to emphasize the relation of these
problems to those of school mathematics. In this way I hope to make it easier
for you to acquire that ability which I look upon as the real goal of your
academic study: the ability to draw (in ample measure) from the great body
of knowledge taught to you here as a vivid stimulus for your teaching.
(Klein, 1908, tr. Schubring, pp. 1--2)
\end{quote}

Thus, this example exemplifies Klein's approach to connecting elementary
and higher mathematics. Yet Klein stops short of advocating changes to
the school curriculum: in schools, students should see only Napier's
treatment, that logarithms are the inverse of exponentiation. This
raises the question of what vertical treatment achieves.

\subsection*{Description of the text}

Klein examines logarithms in increasingly abstract settings, showing
that the school definition is part of a larger structure. He moves from
integers to reals to complex numbers, noting that the latter are needed
to handle negative bases.

Over the positive reals, the logarithm is the inverse of the exponential
function (\(b > 0\)):

\[
\begin{aligned}
y &= b^x, \\
x &= \log_b y.
\end{aligned}
\]

Since exponents satisfy \(b^{c} b^{c'} = b^{c+c'}\), both exponentials
and logarithms connect addition and multiplication fundamentally.

For rational arguments \(y = \frac{m}{n}\), a distinction arises: if
\(n\) is odd, roots exist for any \(b\); if \(n\) is even, roots exist only
for positive \(b\). Hence the curves \(y = \log x\) and \(y = \log(-x)\)
have different densities; the positive half-plane has twice as many points
as the negative.

\begin{figure}[htbp]
\centering
\includegraphics[width=7cm]{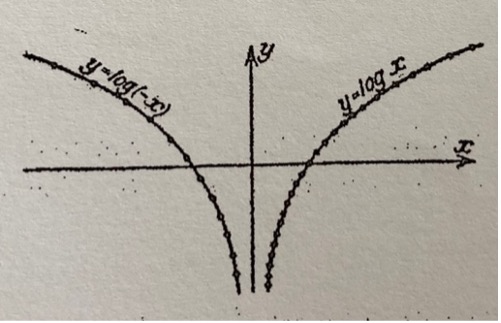}
\caption{Logarithmic and Exponential Functions (EMHS, Schubring p.155)}
\label{fig:logarithm}
\end{figure}

The continuity of both curves involves advanced function theory (complex
analysis), so Klein restricts to \(b > 0\), avoiding conceptual issues
for students. He criticized table-based calculation as ``despicable
utilitarianism'' (p. 156).

Klein traces the historical development: Stifel (1544) produced the first
limited base-2 table; Napier and Bürgi developed the difference equation
underlying logarithms, related to the curve \(1/x\). Following developments
in infinite series (Gauss, Cauchy), he concludes:

\begin{quote}
\small
The introduction of the logarithm by means of the quadrature of the hyperbola
is equal in rigour to any other method, whereas it surpasses all others
in simplicity and clarity. (p. 166)
\end{quote}

Today, the school definition of logarithm is the area under \(1/x\).

\begin{figure}[htbp]
\centering
\includegraphics[width=4cm]{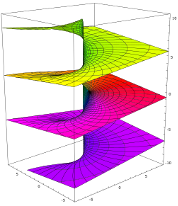}
\caption{Riemann surface of \(\log z\), from Wikipedia}
\label{fig:riemann}
\end{figure}

Klein also discusses the logarithm \(\log z\) as \emph{uniformizing}: a
multivalued function made single-valued and analytic. In the complex plane,
\(\log z\) is multivalued, jumping by \(2\pi i\) each loop around the origin\footnote{
\(\log z = \ln r + i(\theta + 2\pi k), k \in \mathbb{Z}\) for \(z = r e^{i\theta}\).}.
This creates infinitely many branches forming a spiral (Figure~\ref{fig:riemann}),
where the logarithm is well-defined\footnote{Alternatively, one can cut the
plane along the negative real axis and select a branch with angle between
0 and \(2\pi\), flattening one turn of the spiral into a single sheet.}.

Klein’s ``big picture'' shows: over the reals, logarithms simplify
multiplication; over the complex numbers, they unwind the helix for
uniform calculation; in more advanced contexts, elliptic functions
act similarly on toroidal surfaces.

The essential question remains: what is a logarithm fundamentally? Is
its elementary form the same for school students and mathematicians?
Klein raises these questions but does not resolve them.

\subsection*{What this example shows}

Klein’s development follows a vertical trajectory: taking an elementary
concept to a higher level gives a deeper ontological sense, but does not
resolve how the concept should first be presented in schools.

For real numbers, logarithms convert multiplication to addition; for
complex numbers, they straighten the helical structure. Other perspectives,
like power series, exist but are not discussed. The logarithm concept
bridges algebraic and geometric representations.

Klein shows a unifying tendency: logarithms connect discrete arithmetic
progressions and continuous geometric progressions. Yet he accepts the
school-level presentation without introducing the higher viewpoint.

\subsection*{Points to highlight}

\begin{itemize}
    \item Vertical presentation of the concept of logarithm
    \item Historical trajectory showing a topic evolving from calculation to
          research-level applications
    \item Does not solve the double discontinuity problem from the student’s
          perspective
\end{itemize}

%%example 2
\section*{Example 2: Lattice Point Geometry}

The second example we will consider comes from a book called \emph{Roots
to Research: A Vertical Development of Mathematical Problems} (RR) by
Judith and Paul Sally. As the title implies, the treatment is vertical,
but unlike the treatment of logarithms—which takes a single concept and
develops it in increasing layers of complexity—here we look at a variety
of theorems in different fields at different levels of difficulty.

The book has five chapters, each of which focuses on a mathematical
problem that is accessible to young learners and also rich enough to be
explored in research mathematics\footnote{The other topics in this book
  (with some of the advanced developments in parentheses) are the Four
  Numbers problem (Linear Algebra, k-numbers Game), the Rational Right
  Triangles and Congruent Number problem (Elliptic Curves), Rational
  Approximation (Liouville and Thue-Siegel-Roth Theorem), Dissection
  (Hilbert's Third Problem, Banach-Tarski Paradox, Borsuk's Problem).}.
The book is written in the style of a textbook, with theorems and proofs
for each topic, and exercises for students to work on independently. In
the introduction, the authors provide a motivation for the book:

\begin{quote}
\small
Throughout our careers, certain contemporary mathematical problems have
caught our interest because their origins lie in mathematics covered in
the elementary school curriculum, and their development can be traced
through high school, college, and university level mathematics. This
book is intended to provide a source for the mathematics (from beginning
to advanced) needed to understand the emergence and evolution of five of
these problems.
\end{quote}

Because of the structure of the book, all the material in a chapter is
closely related, and theorems are often cross-referenced across chapters
(e.g., the Dirichlet Approximation Theorem appears in two different
chapters). There is occasionally historical context, but the text is
largely developed in the format Lemma, Theorem, Corollary, as in
undergraduate or graduate texts, with historical comments relegated to
footnotes. While the authors intend the book for a ``serious'' high
school student or teachers from high school to graduate level, the
terse style suggests it may be best suited for readers with
mathematical sophistication.

\subsection*{Description of the text}

The chapter begins with a definition of a plane integer lattice—a grid
with horizontal and vertical distances of one unit between adjacent
points. The first sections explore theorems about integer lattices that
do not exceed high school mathematics. For instance, the slope of a line
through two lattice points is rational; the square of the length of a
lattice line segment is an integer; an equilateral triangle cannot be a
lattice polygon.

This leads to Pick's theorem, discussed in three sections, which states
that the area \(A\) of any lattice polygon is given by
\[
A = I + \frac{B}{2} - 1
\]
where \(B\) is the number of boundary points and \(I\) is the number of
interior points. This result is often given to elementary students to
verify for specific cases, as it is initially surprising. RR develops
and proves Pick's theorem in several ways\footnote{See Raman and Ohman,
  2011 for a discussion of different proofs of Pick's theorem.}. For
example, one proof is motivated using Euler's formula for planar graphs,
\(v - e + f = 2\), where \(v\), \(e\), and \(f\) are the numbers of
vertices, edges, and faces, respectively.

The book continues with theorems in lattice geometry, including Farey
sequences, the Dirichlet Approximation Theorem\footnote{The Dirichlet
Approximation Theorem, which is a property of Farey sequences, is not to
be confused with the Dirichlet Prime Number Theorem. The approximation
theorem states: Let \(\alpha\) be a real number and \(n\) a positive
integer. Then there exists a rational number \(\frac{p}{q}\) with \(0 <
q \le n\) such that
\[
\left|\alpha - \frac{p}{q}\right| \le \frac{1}{(n+1)q}.
\]},
lattice points on a circle, integer points in bounded convex regions of
the plane, and Minkowski's theorem. For comparison with the other
examples, we focus here on Farey sequences. These connect lattice
geometry to number theory and deep ideas such as the Riemann Zeta
function. While these topics are advanced, the theorems can be posed and
explored without proof at younger ages. For example, the Polya Orchard
Problem appears as an exercise: ``How thick must the trunks of the trees
in a regularly spaced circular orchard grow so that one could not see
out, standing at the center?''\footnote{See Allen, 1986 for a discussion
  of this problem.}

Farey sequences are defined geometrically: they are the points on a unit
grid ``visible'' from the origin (meaning a line through the origin and
the point contains no other lattice point).

\begin{figure}[htbp]
\centering
\includegraphics[width=5cm]{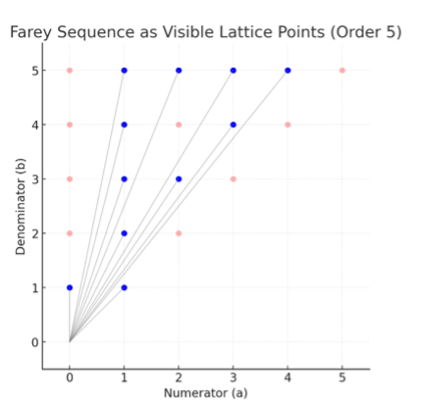}
\caption{Farey sequences as Visible Lattice Points}
\label{fig:farey}
\end{figure}

A point is in a Farey sequence if and only if \(\gcd(a,b) = 1\). For
example, the Farey sequence on a \(5 \times 5\) grid has eleven terms,
which satisfy:

\begin{enumerate}
    \item For consecutive terms \(\frac{a}{b}\) and \(\frac{a'}{b'}\),
          \(a'b - ab' = 1\).
    \item For three consecutive terms \(\frac{a}{b}, \frac{a'}{b'},
          \frac{a''}{b''}\), the mediant property holds:
          \(\frac{a'}{b'} = \frac{a+a''}{b+b''}\).
\end{enumerate}

These facts require only basic linear algebra.

A deeper question explored in the book is the proportion of grid points
that are Farey numbers, equivalent to asking whether two random integers
are coprime. The answer involves the Riemann Zeta function, and the
authors show a clever trick Euler used to calculate it.

\subsection*{What this example shows}

As this example demonstrates, RR contains material spanning elementary to
advanced mathematics. It is satisfying to see how a topic develops from
simple ideas to research-level results. However, the book is not
intended for elementary teachers or students; even the simpler concepts
require some mathematical sophistication to fully understand. The
Corollary-Lemma-Theorem presentation also requires work to unpack.

This is a strong example of vertical development—not just a single
concept, but a small area of mathematics. Like the logarithms example,
the full development is only fully appreciated by someone with
mathematical training. An elementary student or teacher would need
considerable effort to follow it.

\subsection*{Points to highlight}

\begin{itemize}
    \item Vertical presentation of lattice geometry
    \item Begins with school-level problems but develops them primarily through higher-level concepts
    \item Does not resolve the double discontinuity, as mathematical depth is achieved through vertical ascent rather than practice
\end{itemize}

%%example 3

\section*{Example 3: Box Problem}

The third example comes from the textbook \emph{Mathematics for High School
Teachers: An Advanced Perspective} (Usiskin et al.). Throughout this book
we find so-called ``Extended Analyses'' (EA), which start with a routine
high school task and develop it, through careful and systematic inquiry,
into deeper mathematical insights\footnote{For other examples, see pages
5, 76, 125, 173, 197, 375, 440, 492 from MHST. Also see Stanley and
Raman, 2007 for a discussion of Extended Analysis tasks.}. Unlike the
previous examples, these EA problems are horizontal—they do real
mathematics entirely within the bounds of high school mathematics.

The MHST book is intended for pre-service or in-service teachers, or
others who want to study high school mathematics from what the authors
call an ``advanced perspective.'' The chapter explains this term through
three types of extended analyses: (i) concept analysis, illustrated
through parallelism, (ii) problem analysis, illustrated with a task about
averages, and (iii) mathematical connections, illustrated by exploring
the relationship between multiplication and addition, including
exponential functions.

As Stanley notes in his teaching guide\footnote{Dick Stanley, 2003.}:

\begin{quote}
\small
This approach takes advantage of the fact that, implicit in standard
problems of the high school curriculum, there is an undercurrent of
deeper mathematics that is seldom brought out for teachers. Seeing how
to carry out deeper analyses of simple problems gives teachers a better
appreciation of the mathematical substance of the high school
curriculum, and at the same time helps them develop and use powerful
mathematical ways of thinking.
\end{quote}

We will set aside the distinction between deeper and higher mathematics
for the moment and focus on the concrete example from MHST.

\subsection*{Description of the text}

The ``Box Problem'' (section 3.3.3 in MHST), begins with a classic high school
task:

\begin{quote}
\small
Find the maximum volume \(V\) of an open box constructed from an
8-inch by 11-inch cardboard rectangle by cutting squares from the
four corners and folding up the sides.
\end{quote}

Typically, this problem is used for practicing max-min techniques:
express the volume \(V\) as a function of the cut-out side length \(x\),
differentiate, set the derivative equal to zero, solve the quadratic, and
identify the maximum. For this rectangle, the volume function is
\[
V(x) = x(11 - 2x)(8 - 2x),
\]
and the solution is approximately \(x \approx 1.525\) inches, giving
\(V(x) \approx 60\) cubic inches.

This is where a typical high school task stops, but the MHST goes further by suggesting additional analysis:

\begin{quote}
\small
1. Generalize the rectangle sides to \(L\) and \(W\).  

2. Set \(W = 1\) inch and express the volume as a function of \(L\) and
the cut-out length \(x\).  

3. Differentiate to find the maximized volume as a function of \(x\) and
\(L\).  

4. Notice that for a square, the cut-out is \(\frac{1}{6}\) inch and the
base area equals the wall area.  

5. Show that the base area = wall area condition holds for shapes other
than a rectangle.  

6. Show that as \(L\) increases, the limiting cut-out length is
\(x = \frac{1}{4}\) inch.
\end{quote}

These are all prompts during which students can explore and try to find the answers.  Ideally these prompts will come from the students themselves, and the role of the teacher is to keep them on track to find a generalized version of the solution.

Note that the box problem is parallel to another famous calculus problem, the ``pen problem,'' where one maximizes the
area of a rectangular pen against a wall for a fixed material length. The
box problem is the a three-dimensional analog of this two-dimensional
task.

The problem can be explored both algebraically and geometrically. Students
can cut cardboard to build boxes, gaining intuition about which shapes
maximize volume. The key insight is that the maximum occurs when the area
of the base equals the area of the sides. This holds for both 2-D and
3-D cases, and can be proven algebraically by showing that any deviation
from equality reduces the volume.

\subsection*{What this example shows}

This is a high school-level task that leads to genuine mathematical
results. It reflects mathematical practices such as conjecturing,
exploring special cases, analyzing extremes, and seeking structural
explanations, expressed both geometrically and algebraically. The
analysis remains entirely within school-level mathematics, including
simple calculus.

The main difference from the first two examples is that the content level
does not increase. The focus is on methods of inquiry—mathematical
practice itself—rather than moving to higher-level mathematics. Stanley
describes this approach as developing mathematical sophistication with
simple tools:

\begin{quote}
\small
These sessions illustrate an approach to the mathematical preparation
of high school teachers that stays very close to high school content,
but treats it in a mathematically deep and sophisticated way. The intent
is to help teachers focus the kind of mathematical maturity that their
undergraduate courses make possible on the actual content and problems of
high school mathematics. Thus it treats high school mathematics from an
``advanced standpoint.'' It does so through developing the idea of
sophisticated use of simple tools, rather than appealing to higher
levels of rigor or formality.
\end{quote}

We can view this as a separate dimension of development: the first two
examples emphasize vertical content growth, while the Box Problem focuses
on horizontal development of mathematical practices.

\vspace{1em}
Finally, it is important to note that extended analyses differ from
standard ``problem-solving'' or ``teaching for conceptual understanding''
tasks. They do not have a fixed pedagogical goal or time frame. Instead,
they encourage curiosity and mathematical exploration—the very practices
of mathematicians—and can be applied at all levels of the curriculum.

\subsection*{Points to highlight}

\begin{itemize}
\item Horizontal presentation of a school mathematics task
\item Shows richness of high school content, with real results possible
\item Offers a possible solution to the ``double discontinuity'' because
the content level does not exceed calculus, though a teacher must guide
students through the inquiry
\end{itemize}

%%discussion

\section*{Discussion}

The problem of double discontinuity is challenging as long as  we think along one dimension. If ``advanced'' means ``higher'' it is hard to convey full meaning of a concept without understanding the more general context.  As we saw
with the logarithm example, the interesting results emerge only when
considering negative bases (\(b < 0\)), which a typical high school
course would not cover. In the context of lattice geometry, which involves genuine 
mathematical exploration, that exploration also relies on 
relatively advanced concepts, making it out of reach for high school students. The ``extended analysis'' example is quite different from the other
two, because the mathematical richness remains within the
plane of high school mathematics. Using only concepts they know, the mathematics still advances, but it does so via the standard tools of mathematical practice. 

\begin{figure}[htbp]
\centering
\includegraphics[width=5cm]{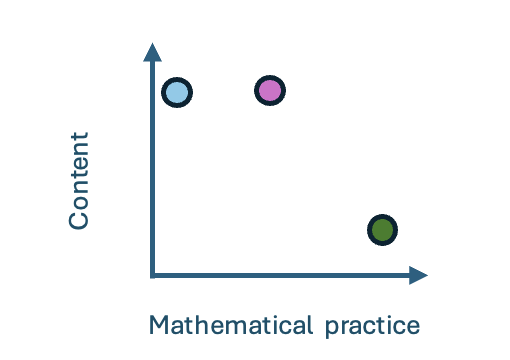}
\caption{Two dimensions of mathematical tasks}
\label{fig:dimensions}
\end{figure}

Comparing these three examples we see that
mathematical richness can proceed along two dimensions. In the first two
cases (vertical), the mathematics develops through content. In the
third case (horizontal), mathematics develops through maturity as
practitioners. This difference can be illustrated with perpendicular
axes, as in Figure~\ref{fig:dimensions}.

Using this diagram, teaching approaches can be represented as a linear
combination of the two axes\footnote{A third dimension involving aesthetics
could be considered, but we leave that for another paper.}. The
logarithms example (light blue) is high on the content axis (consistent
with Kilpatrick and Schubring's insistence on ``higher'' rather than
``advanced''). The lattice geometry example (purple), which involves
discovery, conjecture, and generalization, is fairly high on content
and moderately high on methods. The box problem (green) is very high
on methods, with little change on the content axis.

This diagram reframes the double discontinuity problem, not as a
simple disconnect, but rather a misalignment of dimensions. Focusing only on the vertical
axis makes the problem seem difficult. Allowing for a horizontal
dimension opens a whole new space for solutions.

So what exactly does it mean to present elementary mathematics from
a ``higher standpoint''? Klein himself was not precise about this term.
According to Allmendinger (2016), Klein used it intuitively. Schubring 
(2016) emphasizes its vertical nature, and the image of the teacher being \emph{above} the
subject makes the imagery entirely clear.

However, even in Klein's early writings, we find hints that Klein might have had in mind
some sort of pedagogy beyond the strict vertical development:

\begin{quote}
\small
Instead of developing a proper feeling for mathematical operations, or
promoting a lively, intuitive grasp of geometry, the class-time is spent
learning mindless formalities or practicing pretty trivialities that
exhibit no underlying principle. (Klein, \emph{Antrittsrede}, Rowe tr.,
1985, p.~139)\footnote{%
``Statt den eigentlichen Sinn mathematischer Operationen zu entwickeln,
statt in der Geometrie das lebendige Anschauungsvermögen auszubilden,
wird die Zeit zur Erlernung eines geistlosen Formalismus
oder zur Übung in principlosen Kunststücken verwandt.'' 
(Original, Rowe, 1985, p.~134)
}
\end{quote}

By ``lively''  Klein might have referred to developing a deep understanding of concepts, with (as he puts it) their general principles.  However, one might also read this quote to imply that Klein imagined some sort of real mathematical practice akin to doing research.  Could it be that he did not think about the possibility of doing real, engaging mathematics at the school level?

\subsection*{Coda}
In a sense, nothing in this paper is new. Rather, it is a reframing of what we already know. In modern
mathematics classrooms, particularly at the school level, we are far
more adept at encouraging mathematical practices than in Klein’s day.
Yet the gap between high school and university mathematics persists,
and will continue to exist as long as the curriculum is viewed
primarily in vertical terms. What we see here is that the term
``advanced'' need not be epistemic—that is, tied solely to the
development of concepts along a mathematical trajectory. The term
``advanced'' can also apply to practice, as it does in Stanley’s use.
We learn mathematics not only by appreciating it in its generality,
but also by doing it. ``Advanced'' need not necessarily mean ``higher''.

Klein’s original vision of a school teacher being \emph{above} the
subject does not rule out the possibility that students, too, might be
immersed in genuine mathematical activity. What we have seen here is
that such immersion need not come at the end of a long, hierarchical
ladder of learning. The joy and satisfaction that mathematicians feel
in doing their work can be experienced without leaving the plane of
school mathematics.

%% Appendix

\newpage

\section*{Appendix: Extended Analysis of the Box Problem}

\appendix

A typical box problem reads as follows:

\begin{description}
    \item[Problem 1.] What is the maximum volume $V$ of an open box constructed from a given cardboard rectangle $8'' \times 11''$ by cutting squares out of the corners and folding up the sides?
    
    \item[] This is closely related to the following problem:
    
    \item[Problem 2.] What are the proportions of an open box with a square bottom that maximizes the volume $V$, given that the total area of the sides and bottom of the box is fixed at $88$ square inches?
\end{description}

\noindent
(Note that in Problem 2 the limiting condition is the surface area of the open box, whereas in Problem 1 the limiting condition is the rectangle the box is made from. Problem 2 is arguably more natural since it directly involves the surface area, whereas Problem 1 discards some of the original rectangle.)

There is an extended analysis of Problem 1 in the text \emph{Mathematics for High School Teachers: An Advanced Perspective}. The analysis makes use of the derivative, but the algebra is straightforward, with the side $x$ of the squares cut out from the corners as the unknown. Problem 2, however, involves more variables and a slightly more involved algebraic approach.

\subsection*{Solution to Problem 2}

Let $L$ be the side length of the square bottom and $H$ the height of the sides. The total area of the bottom and side walls is given as $A$ (for the original problem, $A = 88$ in$^2$). The task is to find $L$ and $H$ in terms of $A$ that maximize the volume.

\subsubsection*{Step 1: Express Area in Terms of Unknowns}

\begin{equation}
A = L^2 + 4HL
\end{equation}

\subsubsection*{Step 2: Express Volume in Terms of Unknowns}

\begin{equation}
V = L^2 H
\end{equation}

\subsubsection*{Step 3: Express Volume in Terms of a Single Unknown}

From equation (1), solve for $H$:

\begin{equation}
H = \frac{A - L^2}{4L}
\end{equation}

Substitute into (2):

\begin{align}
V &= L^2 H \\
  &= L^2 \cdot \frac{A - L^2}{4L} \\
  &= L \cdot \frac{A - L^2}{4} \\
  &= \frac{1}{4} (A L - L^3)
\end{align}

\subsubsection*{Step 4: Maximizing Volume}

Differentiate $V$ with respect to $L$:

\begin{equation}
\frac{dV}{dL} = \frac{1}{4} (A - 3L^2)
\end{equation}

Set $\frac{dV}{dL} = 0$:

\begin{equation}
A = 3 L^2
\end{equation}

\noindent
Solve for $L$:

\begin{equation}
L = \sqrt{\frac{A}{3}}
\end{equation}

\subsubsection*{Step 5: Solve for $H$}

From (3):

\begin{equation}
H = \frac{A - L^2}{4L} = \frac{2L^2}{4L} = \frac{L}{2}
\end{equation}

\noindent
Using (7):

\begin{equation}
H = \frac{\sqrt{A/3}}{2}
\end{equation}

\subsubsection*{Step 6: Interpretation}

\begin{equation}
H = \frac{L}{2}
\end{equation}

\noindent
Thus, the open box that maximizes the volume for a given total area $A$ has a square bottom with side length $L$ and wall height $H$ exactly half of $L$. Geometrically, the box is in the shape of half a cube.  

This result mirrors the two-dimensional analog: a rectangular pen using a wall as one side has maximum area when it forms half a square (for a fixed perimeter).  

(See Figure 31, page 127 of MHST: AAP.)

\vspace{1em}
\noindent
\emph{Appendix written by Dick Stanley}

%%References

\section*{References}

\subsection*{Primary sources}

\begin{list}{}%
  {\leftmargin=3pt
   \itemindent=0pt
   \labelwidth=0pt
   \labelsep=0pt
   \parsep=0.2\baselineskip}

\item Klein, F. (1908/2016). \textit{Elementary mathematics from a higher standpoint: Volume I: Arithmetic, algebra, analysis} (G.~Schubring, Trans.). Springer.

\item Klein, F. (1911/2016). \textit{Elementary mathematics from a higher standpoint: Volume II: Geometry} (G.~Schubring, Trans.). Springer.

\item Klein, F. (1924/2016). \textit{Elementary mathematics from a higher standpoint: Volume III: Precision mathematics and approximation mathematics} (M.~Menghini \& A.~Baccaglini-Frank, Trans.; G.~Schubring, Ed.). Springer.

\item Sally, J.~D., \& Sally, P.~J., Jr. (2007). \textit{Roots to Research: A Vertical Development of Mathematical Problems}. American Mathematical Society.

\item Usiskin, Z., Peressini, A.~L., Marchisotto, E., \& Stanley, D. (2003). \textit{Mathematics for high school teachers: An advanced perspective}. Pearson Education.

\end{list}

\subsection*{Additional sources}

\begin{list}{}%
  {\leftmargin=3pt
   \itemindent=0pt
   \labelwidth=0pt
   \labelsep=0pt
   \parsep=0.2\baselineskip}

\item Allmendinger, H. (2019). Examples of Klein's practice: Elementary mathematics from a higher standpoint — Volume I. In H.-G.~Weigand et al.\ (Eds.), \textit{The Legacy of Felix Klein}, 203--213.

\item Boaler, J. (2022). \textit{Mathematical mindsets}. Jossey-Bass (2nd ed.).

\item Cuoco, A., Goldenberg, E.~P., \& Mark, J. (1996). Habits of mind. \textit{Journal of Mathematical Behavior}, 15(4), 375--402.

\item D'Alembert, J.~le Rond (1755). \textit{Éléments des sciences}. \textit{Encyclopédie}, tome~V, 491--497.

\item Fried, M. N., \& Unguru, S. (2001). \textit{The elements of conic sections in the framework of Books I–III}. In \textit{Apollonius of Perga’s Conica: Text, Context, Subtext}. Brill, 222, 56--116.

\item Hardy, G.~H., \& Wright, E. M. (1938). \textit{An introduction to the theory of numbers}. Oxford University Press.

\item Allen, T. (1986). Polya's orchard problem. \textit{The American Mathematical Monthly}, 93(2), 98--104.

\item Kilpatrick, J. (2019). What is and what might be the legacy of Felix Klein? In H.-G. Weigand et al.\ (Eds.), \textit{The Legacy of Felix Klein}, 23--31.

\item Klein, F. (1900). Bericht über die Lage der Schulen und den Unterricht in Mathematik. \textit{Preußisches Kultusministerium}.

\item Klein, F., \& Schimmack, H. (1907). \textit{Der Meraner Lehrplan}. Teubner.

\item NRICH. (n.d.). \textit{NRICH: Mathematics resources}. University of Cambridge. \texttt{https://nrich.maths.org}.

\item Pólya, G. (1945). \textit{How to solve it}. Princeton University Press.

\item Price, E.A., trans. (1911). The ``Meraner Lehrplan''. \textit{The Mathematical Gazette}, 6(95), 179--181.

\item Raman-Sundström, M., \& Öhman, L.-D. (2011). Two beautiful proofs of Pick's theorem. In \textit{Proceedings of CERME 7}, 224--232.

\item Rowe, D. E. (1985). Felix Klein’s “Erlanger Antrittsrede”: A transcription with English translation and commentary. \textit{Historia Mathematica}, 12(2), 123--141.

\item Schubring, G. (2019). Klein's conception of ``elementary mathematics from a higher standpoint.'' In H.-G. Weigand et al.\ (Eds.), \textit{The legacy of Felix Klein}, 33--46. Springer. 

\item Stanley, D. (2003). \textit{Delving deeper}. Unpublished manuscript.

\item Stanley, D., \& Sundström, M. (2007). Extended analyses: There is more to high school math than what you think! Special edition of \textit{Journal of Mathematics Teacher Education}, 10(4--6), 391--397.

\end{list}

\end{document}